\documentclass{article}

\usepackage[letterpaper,margin=0.75in]{geometry}
\usepackage{amsmath,verse}
\usepackage{pdfpages}
\usepackage[shortlabels]{enumitem}
\usepackage{url}
\usepackage[colorlinks=true]{hyperref}
\usepackage[natbibapa,nodoi]{apacite}
\usepackage{tikz}
\usetikzlibrary{calc}
\def\h{0.3}
\def\unit{0.3in}
\tikzset{
  hole/.style={rectangle,fill=black,minimum width=1*\unit,minimum height=2*\h*\unit,inner sep=0,text=white,font=\sffamily}
}

\title{Updating Tom Lehrer's ``The Elements''}
\author{Evan M.\ O'Dorney}

\begin{document}

\maketitle

\begin{abstract}
  We analyze the poetico-comedic criteria by which Tom Lehrer arranged the chemical elements in his 1959 hit ``The Elements,'' a parody to the tune of Gilbert and Sullivan's ``The Modern Major-General'' (1879). We compose a new version to the same tune incorporating all 118 currently known elements and fulfilling, so far as possible, the same artistic aims.
\end{abstract}

Among the voluminous {\oe}uvre of Tom Lehrer (1928--2025), ``The Elements'' stands out as one of his most well-known and best beloved songs, which has entertained generations of schoolchildren since its composition in 1959. It is a parody set to the patter tune ``The Modern Major-General'' from Gilbert and Sullivan's 1879 comic opera \emph{The Pirates of Penzance,} and it has spawned similarly comic settings of other lengthy lists to the same tune, such as the ``Boy Scout Merit Badge Song'' \citep{BoyScout}. For a playful panegyric to Lehrer from a mathematician's perspective, see a recent \emph{Notices} article \citep{Courant_Lehrer}.

``The Elements'' is also one of Lehrer's most \emph{mathematical} songs, not because it is mathematical in subject matter (as are ``New Math,'' ``There’s a Delta for Every Epsilon,'' and others) but because it sets a combinatorial challenge: to permute all the chemical elements into an order fulfilling a number of desiderata, which are not written down but can be inferred from Lehrer's choices, and which interact in nontrivial ways.

At the time ``The Elements'' was written, 102 elements were known (Lehrer included nobelium, whose discovery was still in dispute). This is enough for six $4$-line stanzas of ``The Modern Major-General,'' or two rounds of the full tune. Since then, 16 additional elements have been synthesized, and various performers have attempted to incorporate the new elements into the tune, including Lehrer himself, who composed in 2020 a seventh stanza comprising the 16 new elements. Unfortunately, this seventh stanza lacks the effectiveness of the other six. Lehrer is forced to break the even eighth-note rhythm to squeeze in the longer names of some of the new elements; moreover, a stanza consisting entirely of modern superheavy elements does not quite embody the spirit of ``The Elements,'' which revels in the juxtaposition of the pedestrian and the arcane. Perhaps Lehrer saw the seventh stanza as a stopgap solution, to be revised yet again upon the discovery of yet another element. However, the search for new elements is now in a 16-year hiatus: the latest element, tennessine, was produced in 2009 \citep{ListByDiscovery,Tennessine}, and formally named along with three other recent elements in 2016. Despite numerous attempts around the world, no lab has verifiably synthesized (or mined) a nucleus with more than oganesson's 118 protons, raising the credible possibility that no more cells will appear on the periodic table in our lifetimes, or at all.

Hence the aim of this short note: to compose an updated version of ``The Elements'' arranging the 118 currently known elements in an order consistent with Lehrer's vision and (may I hope) just as entertaining for future generations.

\section{The desiderata}

The lyrics and music to Lehrer's 1959 ``The Elements'' can be accessed for free on Lehrer's website \citep{LehrerTheElements}, which also includes a prominent public-domain declaration relinquishing all rights to his melodies and music. It begins:

\begin{verse}
  There's antimony, arsenic, aluminum, selenium,\\
  And hydrogen and oxygen and nitrogen and rhenium,\\
  And nickel, neodymium, neptunium, germanium,\\
  And iron, americium, ruthenium, uranium, \emph{etc.}
%  
%  Europium, zirconium, lutetium, vanadium,\\
%  And lanthanum and osmium and astatine and radium,\\
%  And gold and protactinium and indium and gallium,\\
%  And iodine and thorium and thulium and thallium.
%  
%  There's yttrium, ytterbium, actinium, rubidium,\\
%  And boron, gadolinium, niobium, iridium,\\
%  And strontium and silicon and silver and samarium,\\
%  And bismuth, bromine, lithium, beryllium, and barium.
%  
%  There's holmium and helium and hafnium and erbium,\\
%  And phosphorus and francium and fluorine and terbium,\\
%  And manganese and mercury, molybdenum, magnesium,\\
%  Dysprosium and scandium and cerium and cesium.
%  
%  And lead, praseodymium and platinum, plutonium,\\
%  Palladium, promethium, potassium, polonium,\\
%  And tantalum, technetium, titanium, tellurium,\\
%  And cadmium and calcium and chromium and curium.
%  
%  There's sulfur, californium and fermium, berkelium,\\
%  And also mendelevium, einsteinium, nobelium,\\
%  And argon, krypton, neon, radon, xenon, zinc and rhodium,\\
%  And chlorine, carbon, cobalt, copper, tungsten, tin and sodium.
%  
\end{verse}

The poem ends with a two-line envoi poking fun at the Boston accent traditionally spoken in the vicinity of Harvard University, where Lehrer once taught:
\begin{verse}
  These are the only ones of which the news has come to Ha'vard,\\
  And there may be many others but they haven't been discavard.
\end{verse}
The envoi is the least mathematical and most readily changeable section of the song; many performers have substituted their own comic couplets here. For the remainder of the article, we will disregard the couplet, focusing on the desiderata that informed the arrangement of the elements.

\begin{enumerate}
  \item \emph{Meter.} Lehrer's work, like Gilbert and Sullivan's original, adheres rigidly to a meter of sixteen syllables per line. The lines can be described as iambic, but it is more correct to point out a three-tiered structure of primary-stressed (\'{}), secondary-stressed (\`{}), and unstressed syllables:
  \begin{verse}
    There's \'antim\`ony, \'arsen\`ic, al\'umin\`um, sel\'eni\`um, \emph{etc.}
  \end{verse}
  Lehrer consistently aligns the primary stress of the element name with a beat of primary stress in the meter (beat 1 or 3 of the measure in the musical score). Occasionally this entails a compromise in the secondary stress (e.g.\ \emph{pras\`eod\'ymi\`um} instead of \emph{pr\`aseod\'ymium}). Elements of one or two syllables, lacking their own secondary stress, can be freely placed on either a primary- or secondary-stressed beat (e.g.\ \emph{\'argon, kr\`ypton, n\'eon, r\`adon,\dots}). (Presumably, Lehrer would have treated an amphibrachic word, such as \emph{paj\'amas,} with the same freedom; but no element has this stress pattern.)
  
  \item \emph{Rhyme.} Lehrer's work, like Gilbert and Sullivan's original, employs \emph{sdrucciola,} that is, triple rhyme wherein the stressed vowel together with the two unstressed syllables following must agree at the ends of each pair of lines (\emph{sel\underline{\'enium},} \emph{rh\underline{\'enium}}). Following the usual tradition in English, there must be a discrepancy in consonants at the onset of the stressed syllable (e.g.\ \emph{ac\underline{t\'inium}/protac\underline{t\'inium}} is avoided). All elements lacking the usual \emph{-ium} ending (e.g.\ \emph{silicon}) have no rhyme in the periodic table and accordingly do not end a line.
  
  \item \emph{Alliteration and other sound similarity.} This is one of the charms of Lehrer's work. Elements beginning with the same sound are agglomerated, regardless of their chemical properties, sometimes spelling a sound in multiple ways (e.g. \emph{\underline{ph}osphorus/\underline{f}rancium}). An especially long chain of seven elements starting with \emph{p} is made possible by the existence of a rhymed pair, \emph{plutonium/polonium.} At other times, successive elements share ending sounds, though not enough to be rhymes (e.g.\ \emph{hydr\underline{ogen}/ox\underline{ygen}/nitr\underline{ogen}}). In the most felicitous cases, there is only one phoneme of discrepancy (\emph{th\underline{u}lium/th\underline{a}llium}, \emph{ce\underline{r}ium/ce\underline{s}ium}).
  
  \item \emph{Continual contrast.} When none of the foregoing similarity properties applies, Lehrer grants us a colorful panoply of elements with unlike sounds:
  \begin{verse}
    Europium, zirconium, lutetium, vanadium,\\
    And lanthanum and osmium and astatine and radium,\\
    And gold and protactinium and indium and gallium, \emph{etc.}
  \end{verse}
  Elements featuring similar sounds such as \emph{m} and \emph{n}, or which nearly rhyme, are diligently kept a safe distance apart. This desideratum weighs more heavily with the addition of the new elements, some of which bear an aurally confusing similarity to existing element names (e.g.\ \emph{nihonium/niobium, hassium/hafnium}).
  
  \item \emph{Placement of the ``and'' connectors.} In English, either ``A, B, and C'' or ``A and B and C'' conveys a complete list of three items, while ``A and B, C'' does not. Some patter songs simply insert ``and'' between every pair of items, such as the Broadway hit about the biblical Joseph's coat of many colors \citep{JosephsCoat}:
  \begin{verse}
    It was red and yellow and green and brown\\
    And scarlet and black and ochre and peach, \emph{etc.}
  \end{verse}
  For Lehrer, however, the meter is overriding: ``and'' is used only when a syllable of space intrudes between two elements. Lehrer tends to prefer \emph{and} at the beginning of lines, where it averts the sense of stoppage, and also before the last element in a stanza: 5 of the 6 stanzas end with \emph{and} and a single element, including all four of the stanzas that precede a piano interlude.
  
  \item \emph{Overall form.} Finally, there are some subtle but effective tactics to drive the form on the scale of the whole piece. Beginning with \emph{a}-alliteration is apt (some listeners will receive the fleeting and mistaken impression of a listing of the elements in alphabetical order). The exciting \emph{th}- and \emph{c}-alliterations are saved for the last lines of stanzas 2 and 5, which follow a fermata in the music. In the sixth and final stanza, there is at first a grouping of the heavier actinoid elements, which add pomp owing to their longer, more consonant-heavy sound:
  \begin{verse}
    There's sulfur, californium and fermium, berkelium,\\
    And also mendelevium, einsteinium, nobelium,
  \end{verse}
  and a final flourish that rattles off two-syllable elements in rapid succession:
  \begin{verse}
    And argon, krypton, neon, radon, xenon, zinc and rhodium,\\
    And chlorine, carbon, cobalt, copper, tungsten, tin and sodium.
  \end{verse}
\end{enumerate}

\section{Considerations in crafting the new version}
When attempting to fit the 16 newest elements into Lehrer's scheme, the main inconvenience is \emph{not} a lack of rhymes; as Lehrer's 2020 stanza shows, there are already two rhymed pairs among the new elements alone (\emph{flerovium/moscovium, livermorium/bohrium}), as well as numerous options linking the old and new elements.

However, fidelity to the meter places a heavy strain on Lehrer's form, specifically finding enough room on secondary-stressed beats (beats 2 and 4 of the musical measure) to accommodate the new elements, all of which have at least one secondary stress and some two, e.g. \emph{r\`utherf\'ordi\`um}. (Here we mark the syllable \emph{-um} as secondary-stressed because the location of the primary stress forces it into this metrical position, even though it would ordinarily be pronounced in an unstressed way.) Indeed, Lehrer's 2020 stanza is too long by three syllables, all of which bear secondary stress. Since we can freely move short elements like \emph{argon} and \emph{krypton} back and forth between primary- and secondary-stressed beats, our first quest is to free up three more stressed beats (primary and secondary together) in the poem.

Two such beats are ready at hand: the filler word ``also'' in the last stanza (``And \underline{also} mendelevium,'' etc.), and the secondary-stressed ``and'' in the third stanza,
\begin{verse}
  And strontium and silicon and silver \underline{and} samarium, \emph{etc.}\\
%   And bismuth, bromine, lithium, beryllium, and barium.
\end{verse}
do not carry an element. It is necessary to free up one more stressed beat. Fortunately, some elements have alternative pronunciations, and in some instances a contracted form with one fewer syllable is possible or common: \emph{fl'orine} (as is nowadays so overwhelmingly standard that singers are liable to trip over Lehrer's three-syllable rendition), \emph{ars'nic,} and a \emph{-shum} pronunciation of elements such as \emph{lutetium, technetium, copernicium,} etc.

% Saving this beat cannot be done without regret, because it prevents the common color-term \emph{silver} from being misapprehended as an adjective, and also because the slight slowing in speed of the succession of \emph{s}-elements sets up, with perfect comic relief, the boisterous barrage of \emph{b}-elements that follows.
% If all else fails, it is feasible to shorten an instrumental interlude by one eighth note and begin a stanza on an extra secondary stress. This is arguably more faithful to the meter than subdividing an eighth note into two sixteenths, as Lehrer recommends in the 2020 edition.
% \end{itemize}
% In the present version, we have opted to eliminate ``also'' and pronounce a few elements in a contracted way. We have retained one secondary-stressed ``and,'' but not in Lehrer's original position.

Once the pronunciation of the elements is fixed, composing the song becomes a kind of tiling problem where the elements must fit in non-overlapping fashion into a grid of stresses (here taller cells denote heavier stress), with any remaining holes to be filled with connectors such as ``there's'' and ``and'':
\[
  \begin{tikzpicture}[x=\unit,y=\unit]
    \foreach \t in {0,4,...,12}{%
      \begin{scope}[shift={(\t,0)}]
        \foreach \hh in {\h, -\h}{%
          \draw (0,\hh) -- (1,\hh) -- (1,3*\hh) -- (2,3*\hh) --
          (2,\hh) -- (3,\hh) -- (3,2*\hh) -- (4,2*\hh) -- (4,\hh);
        }
      \end{scope}
    }
    \foreach \t in {0,1,5,8,12,16}{%
      \draw (\t,-\h) -- (\t,\h);
    }
    \node[hole] at (0.5,0) {\strut\scalebox{.6}{There's}};
    \node  (1) at ( 1.5,0) {\strut an};
    \node  (2) at ( 2.5,0) {\strut ti};
    \node  (3) at ( 3.5,0) {\strut mo};
    \node  (4) at ( 4.5,0) {\strut ny};
    \node  (5) at ( 5.5,0) {\strut ar};
    \node  (6) at ( 6.5,0) {\strut se};
    \node  (7) at ( 7.5,0) {\strut nic};
    \node  (8) at ( 8.5,0) {\strut a};
    \node  (9) at ( 9.5,0) {\strut lu};
    \node (10) at (10.5,0) {\strut mi};
    \node (11) at (11.5,0) {\strut num};
    \node (12) at (12.5,0) {\strut \small \ roent};
    \node (13) at (13.5,0) {\strut ge};
    \node (14) at (14.5,0) {\strut ni};
    \node (15) at (15.5,0) {\strut um};
    
    \path  (1) --  (2) node[midway] {\strut -};
    \path  (2) --  (3) node[midway] {\strut -};
    \path  (3) --  (4) node[midway] {\strut -};
    %    \path  (4) --  (5) node[midway] {\strut -};
    \path  (5) --  (6) node[midway] {\strut -};
    \path  (6) --  (7) node[midway] {\strut -};
    %    \path  (7) --  (8) node[midway] {\strut -};
    \path  (8) --  (9) node[midway] {\strut -};
    \path  (9) -- (10) node[midway] {\strut -};
    \path (10) -- (11) node[midway] {\strut -};
    %    \path (11) -- (12) node[midway] {\strut -};
    \path (12) -- (13) node[midway] {\strut -};
    \path (13) -- (14) node[midway] {\strut -};
    \path (14) -- (15) node[midway] {\strut -};
  \end{tikzpicture}
\]
In jigsaw-like fashion, the longest elements must often be preceded by the shortest ones:
\[
  \begin{tikzpicture}[x=\unit,y=\unit]
    \foreach \t in {0,4,...,12}{%
      \begin{scope}[shift={(\t,0)}]
        \foreach \hh in {\h, -\h}{%
          \draw (0,\hh) -- (1,\hh) -- (1,3*\hh) -- (2,3*\hh) --
          (2,\hh) -- (3,\hh) -- (3,2*\hh) -- (4,2*\hh) -- (4,\hh);
        }
      \end{scope}
    }
    \foreach \t in {0,1,2,8,9,12,16}{%
      \draw (\t,-\h) -- (\t,\h);
    }
    \node[hole] at (0.5,0) {\strut\scalebox{.6}{There's}};
    \node  (1) at ( 1.5,0) {\strut lead};
    \node  (2) at ( 2.5,0) {\strut pra};
    \node  (3) at ( 3.5,0) {\strut se};
    \node  (4) at ( 4.5,0) {\strut o};
    \node  (5) at ( 5.5,0) {\strut dym};
    \node  (6) at ( 6.5,0) {\strut i};
    \node  (7) at ( 7.5,0) {\strut um};
    \node[hole]  (8) at ( 8.5,0) {\strut and};
    \node  (9) at ( 9.5,0) {\strut plat};
    \node (10) at (10.5,0) {\strut i};
    \node (11) at (11.5,0) {\strut num};
    \node (12) at (12.5,0) {\strut plu};
    \node (13) at (13.5,0) {\strut to};
    \node (14) at (14.5,0) {\strut ni};
    \node (15) at (15.5,0) {\strut um};
    
    %\path  (1) --  (2) node[midway] {\strut -};
    \path  (2) --  (3) node[midway] {\strut -};
    \path  (3) --  (4) node[midway] {\strut -};
    \path  (4) --  (5) node[midway] {\strut -};
    \path  (5) --  (6) node[midway] {\strut -};
    \path  (6) --  (7) node[midway] {\strut -};
    %\path  (7) --  (8) node[midway] {\strut -};
    %\path  (8) --  (9) node[midway] {\strut -};
    \path  (9) -- (10) node[midway] {\strut -};
    \path (10) -- (11) node[midway] {\strut -};
    %\path (11) -- (12) node[midway] {\strut -};
    \path (12) -- (13) node[midway] {\strut -};
    \path (13) -- (14) node[midway] {\strut -};
    \path (14) -- (15) node[midway] {\strut -};
  \end{tikzpicture}
\]
One of the new elements to be added has a unique stress pattern: \emph{\`ogan\'esson} (as per \citealt{MWOganesson}; the pronunciation \emph{og\'aness\`on} used by some performers is not so widely accepted). Its unusually late primary stress requires it to be preceded by a short element and followed by a long one (or another short one):
\[
\begin{tikzpicture}[x=\unit,y=\unit]
  \foreach \t in {0,4,...,12}{%
    \begin{scope}[shift={(\t,0)}]
      \foreach \hh in {\h, -\h}{%
        \draw (0,\hh) -- (1,\hh) -- (1,3*\hh) -- (2,3*\hh) --
        (2,\hh) -- (3,\hh) -- (3,2*\hh) -- (4,2*\hh) -- (4,\hh);
      }
    \end{scope}
  }
  \foreach \t in {0,1,3,7,12,16}{%
    \draw (\t,-\h) -- (\t,\h);
  }
  \node[hole] at (0.5,0) {\strut{and}};
  \node  (1) at ( 1.5,0) {\strut i};
  \node  (2) at ( 2.5,0) {\strut ron};
  \node  (3) at ( 3.5,0) {\strut o};
  \node  (4) at ( 4.5,0) {\strut ga};
  \node  (5) at ( 5.5,0) {\strut nes};
  \node  (6) at ( 6.5,0) {\strut son};
  \node  (7) at ( 7.5,0) {\strut am};
  \node  (8) at ( 8.5,0) {\strut er};
  \node  (9) at ( 9.5,0) {\strut ic};
  \node (10) at (10.5,0) {\strut i};
  \node (11) at (11.5,0) {\strut um};
  \node (12) at (12.5,0) {\strut zir};
  \node (13) at (13.5,0) {\strut co};
  \node (14) at (14.5,0) {\strut ni};
  \node (15) at (15.5,0) {\strut um};
  
  \path  (1) --  (2) node[midway] {\strut -};
  %\path  (2) --  (3) node[midway] {\strut -};
  \path  (3) --  (4) node[midway] {\strut -};
  \path  (4) --  (5) node[midway] {\strut -};
  \path  (5) --  (6) node[midway] {\strut -};
  %\path  (6) --  (7) node[midway] {\strut -};
  \path  (7) --  (8) node[midway] {\strut -};
  \path  (8) --  (9) node[midway] {\strut -};
  \path  (9) -- (10) node[midway] {\strut -};
  \path (10) -- (11) node[midway] {\strut -};
  %\path (11) -- (12) node[midway] {\strut -};
  \path (12) -- (13) node[midway] {\strut -};
  \path (13) -- (14) node[midway] {\strut -};
  \path (14) -- (15) node[midway] {\strut -};
\end{tikzpicture}
\]
However, unlike an ordinary tiling problem, the most difficult and rewarding part is arranging the elements with the most common tile shapes (especially \emph{o\'oo\`o} and \emph{\'oo\`o}) to satisfy the other desiderata---rhyme, alliteration, and contrast being the most restrictive. 

\section{The music}
Sullivan's tune to ``I Am the Very Model of a Modern Major-General'' has an AB$\text{A}'$ form. The A and $\text{A}'$ sections are related as antecedent and consequent, in the same major key, and could be sung in uninterrupted succession; but Sullivan separates them with a contrasting section B in the parallel minor. In Sullivan's original, the transition from B back to $\text{A}'$ is prolonged by an echo of the Major-General's last line by the chorus in the dominant key; Lehrer omits this chorus in his adaptation.\footnote{Performers have matched the stanzas of ``The Elements'' to the tune in different ways; since all six stanzas have the same meter, many combinations are possible. The present author has a childhood memory of a version where the first two stanzas are in F major, the third stanza ``There's yttrium\dots'' switches to the F-minor B-section, and at ``There's holmium\dots'' there appears an arrangement of Sullivan's chorus in the dominant key (C), stretched out to four lines by vocalizing the instrumental interlude that leads to the next two verses in the major tonic key. The author would be curious if this version can still be found.} The six stanzas of the 1959 version naturally lend themselves to performance as two identical strophes for an overall form of AB$\text{A}'$-AB$\text{A}'$. With the addition of a seventh stanza, the division into equal strophes cannot be maintained (because $7$ is a prime number!). In place of the rather static AB$\text{A}'$-AB$\text{A}'$-$\text{A}'$ that results when the new stanza is tacked onto an already-complete performance, I propose AAB$\text{A}'$-AB$\text{A}'$.

Lehrer transposed the key of the song from F down to C, which is much more singable by the average person, and I have followed him. Soloists can freely choose a key that fits their individual voice; it is even acceptable to modulate upward for the second strophe (by the standards prevailing in Lehrer's time, though not Gilbert and Sullivan's).

Also, most performers including Lehrer himself simplify the melody, replacing the oscillating neighbor notes with extended patter on a single note. This makes it easier to sing faster, enhancing the comic effect by focusing attention on the rapid-fire lyrics. The present melody is a free transcription of Lehrer's 1967 performance \citep{Lehrer_perform}, with input from his imitators to get a complete melodic line (Lehrer freely dances over the boundary where song meets speech, but for the purposes of group singing, we found it desirable to have a complete melodic line).

As was noted in the memorial of Lehrer's life and music recently published in \emph{Notices} \citep{Courant_Lehrer}, Lehrer publicly relinquished all copyrights and encouraged others to adapt, arrange, alter, and parody his music, just as he had done with the music of his predecessors. I offer the present work in the same spirit. The ordering of the elements in my version is not final. There cannot be a mathematical proof that a particular permutation of the elements is funniest: humor will always rest partially on unexpectedness and novelty. Credit will always be due to Lehrer for hatching the marvelous concept of this song, and to Sullivan for the catchy tune, whereas such an advance as my updating can be encoded in so few bits that any attempt to copyright it would inevitably create an ``illegal prime number'' \citep{IllegalPrime}. Indeed, already, ``The Elements'' has reached a status approaching that of a folk song, with performers making alterations which are iteratively passed on to future generations. This is my contribution; here's hoping it will stand the test of time, at least until someone manages to synthesize element 119.

\bibliographystyle{apacite}
\bibliography{Lehrer_paper}

\includepdf[fitpaper=true, pages=-]{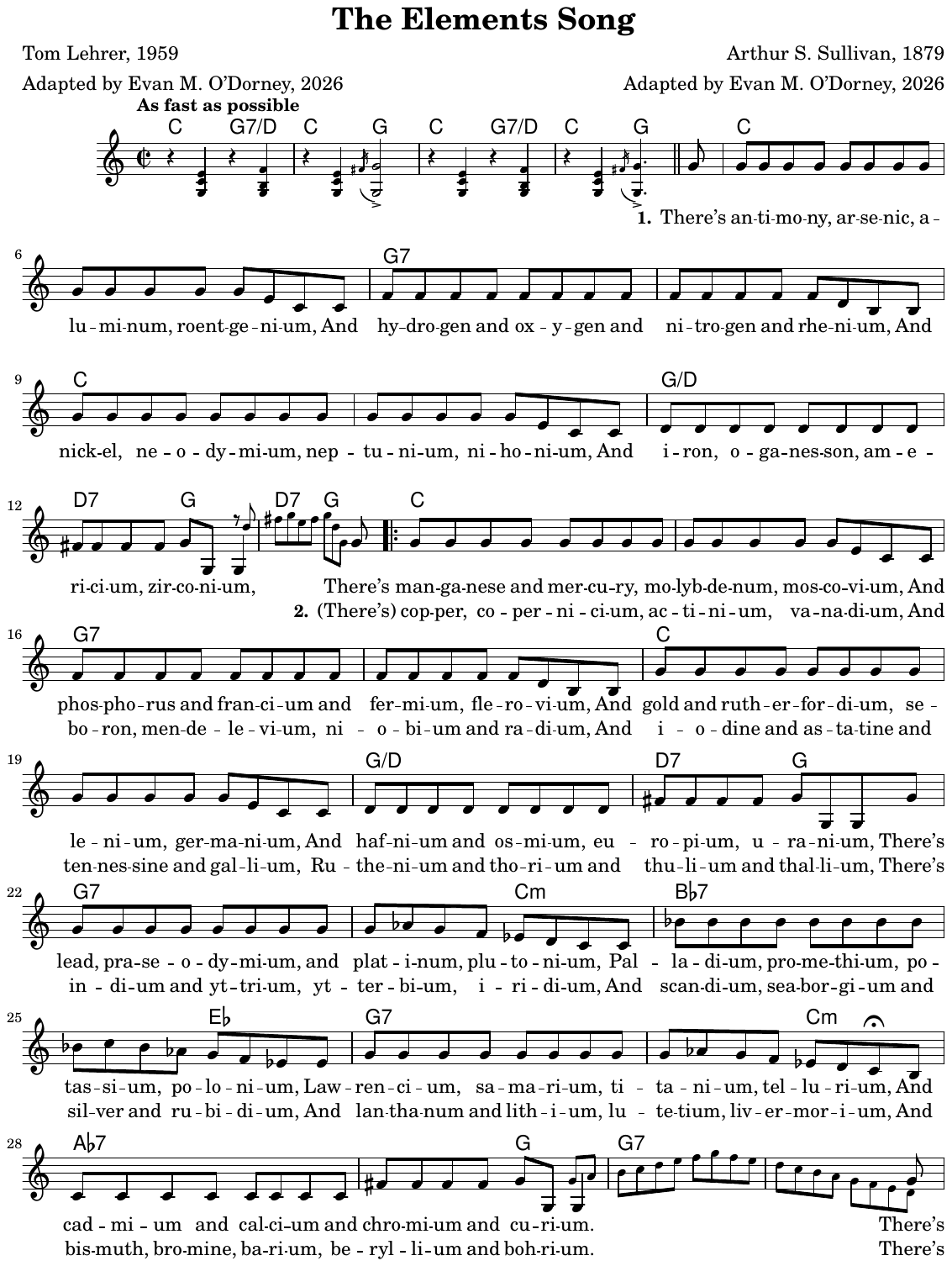}

\end{document}